\documentclass[reqno,12pt]{amsart}
\usepackage{amssymb}
\usepackage{mathrsfs}
\usepackage{amsmath,amssymb}
\usepackage{paralist}

\usepackage{hyperref}
\hypersetup{urlcolor = red, citecolor = blue}
\usepackage[top = 1in,bottom = 1in,left = 1in,right = 1in]{geometry}
\usepackage[sort]{cite}

\allowdisplaybreaks[4]

\newtheorem{theorem}{Theorem}[section]
\newtheorem{corollary}[theorem]{Corollary}

\newtheorem{lemma}[theorem]{Lemma}
\newtheorem{proposition}[theorem]{Proposition}

\theoremstyle{definition}

\newtheorem{remark}[theorem]{Remark}

\numberwithin{equation}{section}

\newcommand{\dd}{\;\mathrm{d}}

\begin{document}

\title[Global solvability and unboundedness in indirect chemotaxis]{Global solvability and unboundedness in a fully parabolic quasilinear chemotaxis model with indirect signal production}  

\author[Mao]{Xuan Mao}
\address{School of Mathematics, Southeast University, Nanjing 211189, P. R. China}
\email{230218181@seu.edu.cn}

\author[Li]{Yuxiang Li}
\address{School of Mathematics, Southeast University, Nanjing 211189, P. R. China}
\email{lieyx@seu.edu.cn}

\thanks{Supported in part by National Natural Science Foundation of China (No. 12271092, No. 11671079) and the Jiangsu Provincial Scientific Research Center of Applied Mathematics (No. BK20233002).}


\begin{abstract}
This paper is concerned with a 
quasilinear chemotaxis model with indirect signal production, $u_t = \nabla\cdot(D(u)\nabla u - S(u)\nabla v)$, $v_t = \Delta v - v + w$ and $w_t = \Delta w - w + u$, posed on a bounded smooth domain $\Omega\subset\mathbb R^n$, subjected to homogenerous Neumann boundary conditions, where nonlinear diffusion $D$ and sensitivity $S$ generalize the prototype $D(s) = (s+1)^{-\alpha}$ and $S(s) = (s+1)^{\beta-1}s$.

Ding and Wang [M.Ding and W.Wang, Discrete Contin. Dyn. Syst. Ser. B, 24 (2019), 4665-4684.] showed that the system possesses a globally bounded classical solution if $\alpha + \beta <\min\{1+2/n,4/n\}$. 
While for the J\"ager-Luckhaus variant of this model, 
namely the second equation replaced by $0=\Delta v - \int_\Omega w/|\Omega| + w$, 
Tao and Winkler [2023, preprint] announced that if $\alpha + \beta > 4/n$ and $\beta>2/n$ for $n\geq3$, 
with radial assumptions, the variant admits occurrence of finite-time blowup. 

We focus on the case $\beta<2/n$, 
and prove that $\beta < 2/n$ for $n\geq2$ is sufficient for global solvability of classical solutions; 
if $\alpha + \beta > 4/n$ for $n\geq4$, 
then radially symmetric initial data with large negative energy enforce blowup happening in finite or infinite time, 
both of which imply that the system allows infinite-time blowup if $\alpha + \beta > 4/n$ and $\beta < 2/n$ for $n\geq 4$.  
\end{abstract}

\maketitle

\section{Introduction}\label{introduce section}

This paper is concerned with global solvability and unboundedness of classical solutions to the fully parabolic quasilinear chemotaxis system with indirect signal production,
\begin{align}
\begin{cases}
	\label{sys: ks isp ppp}
		u_t = \nabla\cdot (D(u) \nabla u) - \nabla \cdot(S(u)\nabla v),&  t>0, x\in\Omega,\\
		v_t =  \Delta v - v + w,&  t>0, x\in\Omega,	\\
		w_t  = \Delta w - w + u, &  t > 0, x\in\Omega, \\
		\partial_\nu u = \partial_\nu v = \partial_\nu w = 0 , & t >0, x\in\partial\Omega,\\
		(u(\cdot, 0), v(\cdot, 0), w(\cdot, 0)) = (u_0, v_0, w_0), & x\in\Omega,
\end{cases}
\end{align}
in a bounded domain $\Omega\subset\mathbb R^n$ with smooth boundary, 
where $\partial_\nu$ denotes the derivative with respect to the outward normal of $\partial\Omega$,
and the initial data $u_0\in C^\vartheta(\overline\Omega)$ for some $\vartheta\in(0,1)$ and 
$u_0, w_0\in W^{1,\infty}(\Omega)$ are nonnegative.
Here, the nonlinear diffusivity $D$ and sensitivity function $S$ generalize the prototype choices 
\begin{equation}
  \label{eq: D S prototype}
  D(s) = K_D(s+1)^{-\alpha}, \quad S(s) = k_S(s+1)^{\beta-1}s, \quad s>0,
\end{equation}
with some $\alpha,\beta\in\mathbb R$ and $K_D,k_S > 0$.


Besides the nature of chemotaxis that describes the oriented movement of cells in response to chemical stimuli \cite{Keller1970},
the system \eqref{sys: ks isp ppp} characterizes on the indirect signal production mechanism \cite{Strohm2013},
that is, the signal $v$ is secreted by the signal producer $w$, rather than the cell $u$ itself,
as the chemotaxis problem with volume-filling effect \cite{Painter2002}
\begin{equation}
  \label{sys: ks pp}
  \begin{cases}
      u_t = \nabla\cdot (D(u) \nabla u) - \nabla \cdot(S(u)\nabla v),&  t>0, x\in\Omega,\\
      v_t =  \Delta v - v + u,&  t>0, x\in\Omega,	\\
      \partial_\nu u = \partial_\nu v = 0 , & t >0, x\in\partial\Omega,\\
      (u(\cdot, 0), v(\cdot, 0)) = (u_0, v_0), & x\in\Omega.
  \end{cases}
\end{equation}
We remark that the system~\eqref{sys: ks isp ppp} is also interpreted as the chemotaxis model with balanced attraction and repulsion \cite{Tao2013,Fujie2017,Guo2018,Li2018}, 
and arises as the main component of a chemotaxis model with phenotype switching \cite{Macfarlane2022,Painter2023,Laurencot2024} as well.  


It has been identified that the chemotaxis models with indirect signal production  
exhibit some distinct features concerning blowup and global existence of solutions, 
compared with the direct ones. 
In the case that the dynamics of signal producers $w$ is rather governed by $w_t = -w + u$,
for instance,  
a novel critical mass phenomenon for infinite-time blowup has been documented by Tao and Winkler~\cite{Tao2017} and Lauren\c{c}ot~\cite{Laurencot2019}, 
and it has been remarked that such a solution blowing up can result in complete infinite-time mass aggregation \cite{Mao2024} (see \cite{Winkler2024,Jin2024} for relevant phenomena).
Moreover, the indirect signal production mechanism has an explosion-suppressed effect 
in the sense that Hu and Tao \cite{Hu2016} showed that 
even weak quadratic logistic dampening can prevent blowup in a three-dimensional chemotaxis-growth model with indirect attraction production,
while Winkler~\cite{Winkler2010} obtained boundedness in the three- and higher-dimensional parabolic-parabolic chemotaxis system with strong logistic source.


We consider the case of fully parabolic chemotaxis models with indirect signal prodcution.
When $\alpha=0$ and $\beta=1$, besides global boundedness in physical spaces $n\leq3$, 
Fujie and Senba~\cite{Fujie2017,Fujie2019} deduced a four-dimensional critical mass phenomenon for the classical solutions of \eqref{sys: ks isp ppp} in the sense that 
\begin{itemize}
  \item  
  under radial assumptions,
  the solution is global and uniform-in-time bounded 
  for nonnegative $(u_0,v_0,w_0)$ with $\int_\Omega u_0 < 64\pi^2$. Without radial assumptions, the same conclusion holds for \eqref{sys: ks isp ppp} subjected to a mixed boundary conditions 
  $\partial_\nu u - u\partial_\nu v = v = w = 0$ on $\partial\Omega$ for $t>0$.
  \item there exist initial data satisfying 
  $\int_\Omega u_0\in(64\pi^2,\infty)\setminus 64\pi^2\mathbb N$ such that 
  the solutions of \eqref{sys: ks isp ppp} with the mixed boundary conditions blow up.
\end{itemize}
Ding and Wang~\cite{Ding2019} found that if $\alpha+\beta<\min\{1+2/n, 4/n\}$,
then for any nonnegative and suitably regular initial datum $(u_0,v_0,w_0)$, 
the problem~\eqref{sys: ks isp ppp} possesses a globally bounded and classical solution.
For the J\"ager-Luckhaus variant of the system~\eqref{sys: ks isp ppp} posed on a ball, 
that is, the second equation in \eqref{sys: ks isp ppp} is replaced by $\Delta v - \int_\Omega w/|\Omega| + w = 0$,
Tao and Winkler~\cite{Tao2023} announced that if $\alpha + \beta > 4/n$ and $\beta > 2/n$ for $n\geq3$,
then there exist radially symmetric initial data such that solutions of the variant blow up in finite time.
It is natural to ask what happens in the case of $\alpha+\beta>4/n$ and $\beta<2/n$.


As motivation, we shall recall two critical lines ($\alpha+\beta=2/n$ and $\beta=0$) appeared in the problem~\eqref{sys: ks pp}.
It is well known that 
\begin{itemize}
  \item if $\alpha + \beta < 2/n$, then the solution of \eqref{sys: ks pp} exists globally and remains bounded \cite{Tao2012}.
  \item if $\alpha + \beta > 2/n$ and $n\geq2$, then for any $m>0$, there exist nonnegative and radially symmetric initial data $(u_0,v_0)$ with $\int_\Omega u_0 = m$ such that solutions to the problem~\eqref{sys: ks pp} posed on a ball, blow up in finite or infinite time \cite{Winkler2010}. 
\end{itemize}
Moreover, built on the method of \cite{Winkler2013}, unboundedness in the case $\alpha+\beta > 2/n$ has been further examined in \cite{Cieslak2012,Cieslak2014,Cieslak2015},
where finite-time blowup is obtained with additional conditions either $\alpha\leq0$ or $\beta\geq1$.
While Winkler \cite{Winkler2019a} showed that whenever $\alpha\in\mathbb R$, $\beta\leq0$ ensures that the system~\eqref{sys: ks pp} has a global classical solution for any suitably regular initial datum, 
and that if $\alpha+\beta>2/n$ ($n\geq2$) and $\beta\leq0$, 
then for any $m>0$, 
there exist initial data $(u_0,v_0)$ with $\int_\Omega u_0 = m$ such that 
the solutions blow up in infinite time 
(see \cite{Lankeit2020} for analogous results for certain parabolic-elliptic version of \eqref{sys: ks pp}).
For the J\"ager-Luckhaus variant of \eqref{sys: ks pp}, 
Winkler and Djie~\cite{Winkler2010a} considered the radially symmetric classical solutions and deduced that 
for any initial data having their mass concentrated sufficiently close to the center of $\Omega$, 
the corresponding solution will undergo a blow-up in finite time, 
provided that $\alpha+\beta > 2/n$ and $\beta>0$. 
Cao and Fuest thus conjectured that $\beta=0$ is the second critical line of \eqref{sys: ks pp}
for finite-time singularity formation in \cite{Cao2024},
where they particularly extended the range of $\alpha$ and $\beta$ for finite-time blowup into 
a subset of $\{(\alpha,\beta): \alpha+\beta>2/n, \beta\in(0,1), \alpha > 0, n\geq2\}$.


The work above suggests us to conjecture that $\beta <2/n$ is sufficient for global solvability of \eqref{sys: ks isp ppp} and infinite-time blowup may happen in the case of $\alpha+\beta>4/n$ and $\beta<2/n$.
To describe our results, 
we assume 
\begin{equation}
  \label{h: D and S}
  \begin{cases}
    D\in C^2([0,\infty)) \text{ satisfies } D > 0 \text{ in } [0,\infty),\quad \text{and that}\\ 
    S\in C^2([0,\infty)) \text{ is positive in } (0,\infty) \text{ and } S(0) = 0,
  \end{cases}
\end{equation}
and that the diffusivity $D$ decays at most algebraically
\begin{equation}
  \label{h: D decay}
  D(s) \geq k_D(1+s)^{-M}, 
  \quad \text{for all } s > 0,
\end{equation}
for an application of a Moser-type iteration, 
where $k_D>0$ and $M\in\mathbb R$ are constants.


Based on a direct $L^p$ estimate for $u$ in case of suitably weak sensitivity strength, 
via maximal Sobolev regularity theories, 
global classical solvability can be asserted by a well-established Moser-type iteration.

\begin{theorem}
  \label{thm: global solvability}
  Let $n\geq2$ and $\Omega\subset\mathbb R^n$ be a bounded domain with smooth boundaries.
  Assume that  $D$ and $S$ comply with \eqref{h: D and S} and \eqref{h: D decay}.
  If there exist $K_S>0$ and 
  \begin{displaymath}
    \beta < \frac2n
  \end{displaymath}
  such that $S(s)\leq K_S(1+s)^{\beta-1}s$ holds for all $s>0$, 
  then for any nonnegative 
  $(u_0,v_0,w_0)\in C^\vartheta(\overline\Omega)\times W^{1,\infty}(\Omega)\times W^{1,\infty}(\Omega)$ 
  with $\vartheta\in(0,1)$, 
  the corresponding classical solution $(u,v,w)$ of the system \eqref{sys: ks isp ppp} exists globally.
\end{theorem}

Concerning possible infinite-time blowup, we extend Winkler's work \cite{Winkler2010} for \eqref{sys: ks pp} to the following 
with considerable modification,
e.g., involving in estimates on a fourth-order linear elliptic equation by certain Hardy-Rellich inequality, 
and initial data constructed with an emphasis on compatible boundary conditions.

\begin{theorem}\label{thm: blowup}
  Let $n\geq4$ and $\Omega \subset\mathbb R^n$ be a ball.
  Assume that  $D$ and $S$ comply with \eqref{h: D and S} and \eqref{h: D decay},
  and suppose that there exist $s_0 > 1$, $\varepsilon\in(0,1)$, $K > 0$ and $k > 0$ such that 
  \begin{equation}
    \label{h: conditions for low bounds of stationary energy}
    \int_{s_0}^s\frac{\sigma D(\sigma)}{S(\sigma)}\dd\sigma 
\leq 
    \begin{cases}
      \frac{Ks}{\ln s},
      & \text{if } n = 4,\\
      \frac{n-4-\varepsilon}{n} \int_{s_0}^s\int_{s_0}^\sigma\frac{D(\tau)}{S(\tau)}\dd\tau\dd \sigma + Ks, 
      & \text{if } n > 4,
    \end{cases}
  \end{equation}
  as well as 
  \begin{equation}
    \label{h: conditions for large negative initial energy}
    \int_{s_0}^s\int_{s_0}^\sigma\frac{D(\tau)}{S(\tau)}\dd\tau\dd \sigma 
    \leq 
    \begin{cases}
      ks(\ln s)^\theta, &  \text{ with some } \theta\in(0,1) \text{ if } n = 4,\\
      ks^{2-\gamma}, &  \text{ with some } \gamma > \frac{4}{n} \text{ if } n > 4,
    \end{cases}
  \end{equation}
  holds for all $s\geq s_0$.
  Then for each $m>0$, there exist radially symmetric initial data $(u_0, v_0, w_0)\in (C^\infty(\overline{\Omega}))^3$ 
  with $\int_\Omega u_0 = m$ such that the corresponding solution $(u,v,w)$ blows up in either finite or infinite time.
\end{theorem}

\begin{remark}
  (i). Inspired by previous work, 
  we detect unboundedness relying on a Lyapunov functional associated to \eqref{sys: ks isp ppp}, 
  which links energy of initial datum with energy of certain steady state 
  via dissipation, for any globally bounded solution.
  Stationary energy has uniformly lower bounds for any steady state with mass of $u$ fixed,
  whereas a family of initial data can be constructed such that their energies has no lower bound,
  under conditions \eqref{h: conditions for low bounds of stationary energy} 
  and \eqref{h: conditions for large negative initial energy}, respectively.
  The conclusion follows by absurdum.

  (ii). This fails to work for \eqref{sys: ks isp ppp} in physical dimensions $n\leq 3$, 
  because initial energy has lower bounds for any reasonably regular initial datum with fixed mass, 
  see Proposition~\ref{prop: mathcal F low bounds in physical space}.

  (iii). The underlying energy method was introduced by Horstmann and Wang~\cite{Horstmann2001} for the classical chemotaxis model, 
  applied by Jin and Wang~\cite{Jin2016} to competing chemotaxis, 
  by Fujie and Jiang~\cite{Fujie2021} to chemotaxis models with density-suppressed motilities, 
  by Lauren\c{c}ot and Stinner~\cite{Laurencot2021} to chemotaxis models with indirect signal production,
  adapted by Winkler~\cite{Winkler2010} to chemotaxis models with volume-filling effect, 
  by Rani and Tyagi~\cite{Rani2024} to chemotaxis-haptotaxis models.
\end{remark}

Applying two theorems above to the system~\eqref{sys: ks isp ppp} with the prototype choices \eqref{eq: D S prototype}, we immediately obtain existence of infinite-time blowup in the case of $\alpha+\beta>4/n$ and $\beta<2/n$.

\begin{corollary}
  \label{coro: infinite-time blowup}
  Let $n\geq4$ and $\Omega\subset\mathbb R^n$ be a ball.
  Suppose that $D$ and $S$ are given by \eqref{eq: D S prototype}.
  If 
  \begin{equation*}
    \alpha + \beta > \frac{4}{n},
  \end{equation*}
  then the conclusion in Theorem~\ref{thm: blowup} holds.  
  If moreover 
  \begin{equation*}
    \beta < \frac2n,
  \end{equation*}
  then the solution $(u,v,w)$ blows up in infinite time.
\end{corollary}

This paper is organized as follows. In Section~\ref{section preliminary},
we introduce some preliminaries, 
including the energy functional and the stationary problem associated with the system~\eqref{sys: ks isp ppp}, 
as basic tools used to detect singularity formation.
Section~\ref{sec: global solvability} is devoted to the proof of Theorem~\ref{thm: blowup}.
We give lower bounds for energy of steady states in Section~\ref{sec: lower bounds for stationary energy}.
We construct a family of initial data with arbitrarily large negative energy 
and thereby complete the proof of Theorem~\ref{thm: blowup} in Section~\ref{sec: initial data with large negative energy}.

\section{Preliminaries. Energy functional and stationary problem}
\label{section preliminary}

We first revisit the results on local existence and uniqueness of the classical solution to the system \eqref{sys: ks isp ppp}.

\begin{lemma}
  \label{le: local existence and uniqueness}
 Let $\Omega \subset \mathbb{R}^n$ be a bounded domain with smooth boundary. 
 Assume that $D, S \in C^2([0, \infty))$ satisfy $D(s)>0$ for $s \geq 0$ and $S(0)=0$. 
 Furthermore, suppose that $u_0 \in C^\vartheta(\overline{\Omega})$ for some $\vartheta\in(0,1)$ and
 $v_0, w_0 \in W^{1, \infty}(\Omega)$ are nonnegative. 
 Then there exist $T_{\max } \in(0, \infty]$ and a unique triplet $(u, v, w)$ of nonnegative functions from $C^0\left(\bar{\Omega} \times\left[0, T_{\max }\right)\right) \cap C^{2,1}\left(\bar{\Omega} \times\left(0, T_{\max }\right)\right)$ solving \eqref{sys: ks isp ppp} classically in $\Omega \times\left(0, T_{\max }\right)$, and if $T_{\max }<\infty$,
  \[
  \limsup _{t \rightarrow T_{\max }}\|u(\cdot, t)\|_{L^{\infty}(\Omega)}=\infty.
  \]
Furthermore,
\begin{equation}
  \label{eq: mass identity}
	\int_{\Omega} u(x, t) \dd {x}=\int_{\Omega} u_{0} \dd x
    \quad\text{for all } t > 0.
\end{equation}
In particular, if $\Omega = B_R := \{x\in\mathbb R^n: |x|<R\}$ for some $R>0$,
and $(u_0, v_0, w_0)$ is a triplet of radially symmetric functions, 
then $u$, $v$ and $w$ are all radially symmetric.
\end{lemma}

\begin{proof}
  The local existence, nonnegativity and the extensibility criterion of classical solutions to the system~\eqref{sys: ks isp ppp} was proven as in \cite[Lemma~3.1]{Ding2019}. 
  The uniqueness can be asserted by Amann's abstract theories \cite[Theorem~14.4]{Amann1993}.
  The mass identity immediately follows from integrating the first equation in \eqref{sys: ks isp ppp}.
  Conservation of radial symmetry is a consequence of uniqueness of solutions and the adequate form of equations in \eqref{sys: ks isp ppp}.
\end{proof}

As observed in \cite{Fujie2017}, the system~\eqref{sys: ks isp ppp} admits a Lyapunov functional,
which plays a key role in detecting blowup. 
Here and henceforth, denote 
\begin{equation}
  G(s) := \int_{s_0}^s\int_{s_0}^\sigma\frac{D(\tau)}{S(\tau)}\dd\tau\dd \sigma, \quad s>0,
\end{equation}
for some $s_0 > 0$. Without loss of generality, we may assume that $G$, as the nature of energy functional, is defined to be positive on $(0,\infty)$.

\begin{lemma}
  \label{le: Lyapunov functional}
  Assume that $D$ and $S$ comply with \eqref{h: D and S}. 
  For the solution $(u,v,w)$ and the maximal time $T_{\max}\in(0, \infty]$ given by Lemma~\ref{le: local existence and uniqueness},
  it holds that 
  \begin{equation}
    \label{eq: Lyapunov functional}
    \frac{\dd}{\dd t}\mathcal F(t) = - \mathcal D(t), \quad t\in(0, T_{\max}),
  \end{equation}
  where 
  \begin{equation}
    \label{eq: mathcal F}
    \mathcal F(t) := \mathcal F(u, v, w) 
    = \int_\Omega G(u) - \int_\Omega uv 
    + \frac12\int_\Omega v_t^2 + \frac{1}{2}\int_\Omega |-\Delta v + v|^2,
  \end{equation}
  and 
  \begin{equation*}
    \mathcal D(t) := \mathcal D(u, v, w) 
    = \int_\Omega S\left\lvert \frac{D}{S}\nabla u - \nabla v\right\rvert^2 
      + 2\int_\Omega |\nabla v_t|^2 + 2\int_\Omega v_t^2.
  \end{equation*}
  Moreover, suppose that $v_0 \in W^{2,2}(\Omega)$,
  with $\partial_\nu v_0 = 0$ on the boundary $\partial\Omega$ in the sense of trace,
  then 
  \begin{equation}
    \label{eq: int Lyapunov functional}
    \mathcal F(u,v,w) - \mathcal{F}_0(u_0,v_0,w_0) = - \int_0^t\mathcal D(\zeta)\dd\zeta, 
    \quad t \in (0, T_{\max}),
  \end{equation}
  with 
  \begin{equation*}
    \mathcal{F}_0(u_0,v_0,w_0) 
    = \int_\Omega G(u_0) - \int_\Omega u_0v_0 
    + \frac12\int_\Omega |\Delta v_0 - v_0 + w_0|^2 
    + \frac{1}{2}\int_\Omega |-\Delta v_0 + v_0|^2.
  \end{equation*}
\end{lemma}

\begin{proof}
  Integrating by parts, we compute 
  \begin{equation}
    \label{eq: diff G - uv}
    \begin{aligned}
    \frac{\dd }{\dd t}\int_\Omega (G(u) - uv) 
    &= \int_\Omega (G'(u) - v)\nabla\cdot(D(u)\nabla u - S(u)\nabla v) - \int_\Omega uv_t\\
    &= - \int_\Omega S\left\lvert \frac{D}{S}\nabla u - \nabla v\right\rvert^2 - \int_\Omega uv_t,
  \end{aligned}
\end{equation}
By the third equation $u = w_t - \Delta w + w$ in \eqref{sys: ks isp ppp}, we have 
\begin{align*}
  \int_\Omega uv_t 
  &= \int_\Omega(w_t - \Delta w + w)v_t = \int_\Omega w_tv_t + \int_\Omega(-\Delta v_t + v_t)w.
\end{align*}
By the second equation $w = v_t - \Delta v + v$ in \eqref{sys: ks isp ppp} and its differential $w_t = v_{tt} - \Delta v_t + v_t$, 
we calculate 
\begin{align*}
  \int_\Omega (-\Delta v_t + v_t)w 
  & = \int_\Omega (-\Delta v_t + v_t)(v_t - \Delta v + v)\\
  & = \int_\Omega |\nabla v_t|^2 
  + \int_\Omega v_t^2 
  + \frac{1}{2}\frac{\dd}{\dd t}\int_\Omega |-\Delta v + v|^2
\end{align*}
and
\begin{align*}
  \int_\Omega w_tv_t 
  & = \int_\Omega (v_{tt} - \Delta v_t + v_t)v_t 
  = \frac12\frac{\dd}{\dd t}\int_\Omega v_t^2 
  + \int_\Omega |\nabla v_t|^2
  + \int_\Omega v_t^2.
\end{align*}
Collecting three identities above, we obtain 
\begin{equation}
  \label{eq: int uvt}
  \begin{aligned}
    \int_\Omega uv_t 
    &=  2\int_\Omega |\nabla v_t|^2 
    + 2\int_\Omega v_t^2 
    + \frac12\frac{\dd}{\dd t}\int_\Omega v_t^2 
    + \frac{1}{2}\frac{\dd}{\dd t}\int_\Omega |-\Delta v + v|^2.
  \end{aligned}
\end{equation}
Substituting \eqref{eq: int uvt} into \eqref{eq: diff G - uv}, we deduce \eqref{eq: Lyapunov functional} as desired.
Thanks to Amann's abstract theories \cite[Corollary~14.5]{Amann1993}, we have 
\[
v\in C([0,T_{\max}); W^{2,2}(\Omega)) \cap C^1([0,T_{\max}); L^2(\Omega)),
\]
which enables us to integrate \eqref{eq: Lyapunov functional} over $(0, t)$ for $t\in(0,T_{\max})$ and get \eqref{eq: int Lyapunov functional} by Newton-Leibniz formula.
\end{proof}

 Let us introduce the stationary problem of \eqref{sys: ks isp ppp} reading as  
\begin{equation}
  \label{sys: stationary system}
  \begin{cases}
    0 = D(u) \nabla u - S(u)\nabla v, &  x\in\Omega,\\
    0 =  \Delta v - v + w, & x\in\Omega,	\\
    0 = \Delta w - w + u, & x\in\Omega, \\
    0 = \partial_\nu v = \partial_\nu w, & x\in\partial\Omega, 
  \end{cases}
\end{equation}
and denote the set of stationary solutions with mass $m$ as 
\begin{equation}
  \label{eq: set of stationary solutions}
  S_m :=\left\{(u,v,w)\in (C^2(\overline{\Omega}))^3: (u,v,w) \text{ solves } \eqref{sys: stationary system}\text{ classically, with } \int_\Omega u = m\right\}.
\end{equation}
For a globally bounded classical solution, 
the key observation is the formula \eqref{eq: F initial - F steady = - int D} 
that the stationary energy is linked to the initial energy via dissipation.
\begin{lemma}
  \label{le: link initial energy to stationary energy}
  Assume that $D$ and $S$ comply with \eqref{h: D and S}.
  Suppose that $(u,v,w)$ is a globally bounded classical solution to the system~\eqref{sys: ks isp ppp}, with $\int_\Omega u_0 = m$, given by Lemma~\ref{le: local existence and uniqueness}.
  Then there exist a sequence $\{t_j\}_{j\in\mathbb N}$ and a triplet $(u_\infty,v_\infty,w_\infty)\in S_m$ such that 
  $t_j\to\infty$ and 
  \begin{equation*}
    (u(\cdot, t_j),v(\cdot,t_j),w(\cdot,t_j))\to(u_\infty,v_\infty,w_\infty)\quad \text{in } \left(C^2(\Omega)\right)^3,
  \end{equation*} 
  as $j\to\infty$, and if moreover $(u_0,v_0,w_0)\in(C^2(\overline{\Omega}))^3$ 
  with $\partial_\nu u_0 = \partial_\nu v_0 = \partial_\nu w_0 = 0$ on $\partial\Omega$, then  
  \begin{equation}
    \label{eq: F initial - F steady = - int D}
    \mathcal{F}(u_\infty,v_\infty,w_\infty) - \mathcal{F}_0(u_0,v_0,w_0) = - \int_0^\infty \mathcal{D}(t)\dd t,
  \end{equation}
  where $\mathcal F$, $\mathcal F_0$ and $\mathcal D$ are given as in Lemma~\ref{le: Lyapunov functional}.
\end{lemma}

\begin{proof}
  Since $(u,v,w)$ is uniform-in-time bounded, the standard parabolic regularity theories entail that 
  there exists $C > 0$ such that 
  \begin{equation*}
    \|(u,v,w)\|_{C^{2+\sigma,1+\sigma/2}(\Omega\times(t,t+1))} < C, \quad t > 1,
  \end{equation*}
  holds for some $\sigma\in(0,1)$. 
  This deduces the boundedness of $F(u,v,w)$ for $t>1$ and thereby $\mathcal D$ is integrable over $(1,\infty)$.
  Besides, by Arzel\`a-Ascoli theorem, 
  there exist a sequence $\{t_j: t_j\geq j, j\in\mathbb N\}$ and a triplet $(\tilde u,\tilde v,\tilde w)\in C^{2,1}(\overline\Omega\times[0,1])$ such that  
  \begin{align*}
    (u_j,v_j,w_j)
    &:= (u(x, t+t_j), v(x, t+t_j), w(x,t+t_j)) \\
    &\to (\tilde u,\tilde v,\tilde w),
    \quad \text{in } \left(C^{2,1}(\Omega\times(0,1))\right)^3,
  \end{align*}
  as $j\to\infty$.
  Abbreviating $\Omega_0 := \Omega\times(0,1)$ 
  and $\{\tilde u > 0\} := \{(x,t)\in\Omega_0: \tilde u(x,t) > 0\}$, 
  which is nonempty and relatively open with respect to $\Omega_0$ due to the mass identity \eqref{eq: mass identity} and continuity of $u$,
  we thus have  
  \begin{align*}
    0 &= \lim_{j\to\infty}\int_{t_j}^{t_j+1}\mathcal D(t)\dd t 
    = \lim_{j\to\infty}\int_{\Omega_0} S(u_j)\left\lvert \frac{D(u_j)}{S(u_j)}\nabla u_j - \nabla v_j\right\rvert^2\\ 
    &\quad + 2\lim_{j\to\infty}\int_{\Omega_0} |\nabla v_{jt}|^2 + 2\lim_{j\to\infty}\int_{\Omega_0} v_{jt}^2\\
    &= I + 2\int_{\Omega_0} |\nabla \tilde v_{t}|^2 + 2\int_{\Omega_0} \tilde v_{t}^2,
  \end{align*}
  with 
  \begin{align*}
    I &:= \lim_{j\to\infty}\int_{\Omega_0} S(u_j)\left\lvert \frac{D(u_j)}{S(u_j)}\nabla u_j - \nabla v_j\right\rvert^2\\
    &\geq \limsup_{j\to\infty}\int_{\{\tilde u>0\}} S(u_j)\left\lvert \frac{D(u_j)}{S(u_j)}\nabla u_j - \nabla v_j\right\rvert^2\\
    &\geq \int_{\{\tilde u>0\}} S(\tilde u)\left\lvert \frac{D(\tilde u)}{S(\tilde u)}\nabla \tilde u - \nabla\tilde v\right\rvert^2
  \end{align*}
  by Fatou lemma. Therefore we get $\tilde v_t = \Delta\tilde v - \tilde v + \tilde w \equiv 0$ over $\Omega_0$, 
  which implies $\tilde w_t = \Delta\tilde w - \tilde w + \tilde u\equiv 0$ over $\Omega_0$,
  as well as 
  \begin{equation}
    \label{eq: D/S grad u = grad v on tilde u > 0}
    \frac{D(\tilde u)}{S(\tilde u)}\nabla \tilde u = \nabla\tilde v \quad \text{over }\{\tilde u > 0\}.
  \end{equation}
  We claim $\{\tilde u>0\} = \Omega_0$ 
  and thereby $(\tilde u,\tilde v,\tilde w)\in S_m$ for each $t\in[0,1]$. 
  Let $\mathcal C$ be a connected component of $\{\tilde u > 0\}$ and denote 
  \begin{equation*}
    f(s) := \int^\varepsilon_s\frac{D(\tau)}{S(\tau)}\dd\tau, \quad s>0,
  \end{equation*}
  with $\varepsilon\in(0,1)$. 
  By \eqref{eq: D/S grad u = grad v on tilde u > 0},
  we have for some $L\in\mathbb R$ such that  
  \begin{displaymath}
    L - \tilde v = f(\tilde u),
  \end{displaymath}
  holds for each $(x,t)\in \mathcal C$.
  Using \eqref{h: D and S}, 
  we can infer $C_\varepsilon>0$ such that 
  \begin{displaymath}
    \frac{D(\tau)}{S(\tau)} \geq \frac{C_\varepsilon}\tau,\quad \tau\in(0,\varepsilon).
  \end{displaymath}  
  Let us suppose by absurdum that $(x_0,t_0)\in\partial\mathcal C\cap\Omega_0$. 
  Then there exists a family $\{(x_i,t_i)\in \mathcal C\}$ such that $(x_i,t_i)\to(x_0,t_0)$.
  Since $u(x_0,t_0)=0$, we may assume without loss of generality that $u(x_i,t_i)\in(0,\varepsilon)$ holds for each $i\in\mathbb N$.
  We thus get 
  \begin{displaymath}
    L - \tilde{v}(x_0,t_0)  
    = L -\lim_{i\to\infty} \tilde{v}(x_i,t_i)  
    = \lim_{i\to\infty} f(\tilde u(x_i,t_i)) 
    \geq \liminf_{i\to\infty}\int_{u(x_i,t_i)}^\varepsilon\frac{C_\varepsilon}{\tau}\dd\tau = +\infty,
  \end{displaymath}
  which is incompatible with the boundedness of $\tilde v$.
  We therefore have $\partial\mathcal C\cap\Omega_0=\emptyset$ and $\mathcal C$ is a relatively closed subset of $\Omega_0$. 
  It follows by connectivity that $\Omega_0 = \mathcal C \subset \{\tilde u>0\}$.
  We thus finish the proof by putting $(u_\infty,v_\infty,w_\infty) = (\tilde u(\cdot,0),\tilde v(\cdot, 0),\tilde w(\cdot ,0))$, 
  applying \eqref{eq: int Lyapunov functional} with $t=t_j$,
  and taking a limit to deduce \eqref{eq: F initial - F steady = - int D}.
\end{proof}

We need the Hardy-Rellich inequality for $n=4$, 
see e.g. \cite[Theorem~3.5 and Corollary~3.2 for the case of radial functions]{Ghoussoub2011} for details, 
\cite{Cazacu2020} for a different proof, 
as well as \cite{Gesztesy2024} and reference therein for general Hardy-Rellich inequalities.

\begin{lemma}
  \label{le: Hardy-Rellich inequality n=4}
  Let $n=4$ and $\Omega = B_R\subset\mathbb R^n$ for some $R>0$.
  Then it holds that
  \begin{equation*}
    \label{eq: Hardy-Rellich inequality n=4}
    \int_\Omega |\Delta u|^2 \geq 3\int_\Omega \frac{|\nabla u|^2}{|x|^2},
  \end{equation*}
  for all $u\in C_0^\infty(\Omega)$,
  and moreover,
  \begin{equation*}
    \label{eq: Hardy-Rellich inequality n=4 for radial functions}
    \int_\Omega |\Delta u|^2 \geq 4\int_\Omega \frac{|\nabla u|^2}{|x|^2},
  \end{equation*}
  for all radial functions $u\in C_0^\infty(\Omega)$.
\end{lemma}

We modify slightly the Hardy-Rellich inequality above to fit our applications.

\begin{lemma}
  \label{le: modified Hardy-Rellich inequality n=4}
  Let $n=4$ and $\Omega = B_R\subset\mathbb R^n$ for some $R>0$. 
  Then 
  \begin{equation}
    \label{eq: modified Hardy-Rellich inequality n=4}
    \int_\Omega \frac{|\nabla u|^2}{|x|^2} 
    \leq \frac{1}{4} \int_\Omega |\Delta u|^2,
  \end{equation}
  holds for all radial functions $u\in C^2(\overline{\Omega})$ with $\partial_\nu u = 0$ over $\partial\Omega$.
\end{lemma}

\begin{proof}
  Let $\tilde{u} := u - u(R)$. Then $\tilde{u} \in H_0^2(\Omega)$. 
  Since $C_0^\infty(\Omega)$ is a dense subset of $H_0^2(\Omega)$, 
  we can find a family $\{u_k\}_{k\geq1}\subset C_0^\infty(\Omega)$ such that 
  \begin{equation*}
    u_k \to \tilde u \quad \text{in } H^2(\Omega),
  \end{equation*}
  as $k\to\infty$.
  Then for each $\delta\in(0, R)$, writing $L^2_\delta := L^2(\Omega\setminus B_\delta)$, we have 
  \begin{align*}
    \left\| \frac{\nabla u}{|x|}\right\|_{L^2_\delta} 
     = \left\| \frac{\nabla\tilde u}{|x|}\right\|_{L^2_\delta} 
    \leq \left\| \frac{\nabla(\tilde u - u_k)}{|x|}\right\|_{L^2_\delta} 
    + \left\| \frac{\nabla u_k}{|x|}\right\|_{L^2_\delta}  
    \leq \frac1\delta\left\| \nabla(\tilde u - u_k)\right\|_{L^2(\Omega)} 
    + \left\| \frac{\nabla u_k}{|x|}\right\|_{L^2(\Omega)}.
  \end{align*}
  Applying Lemma~\ref{le: Hardy-Rellich inequality n=4} to the last term of the inequality above, 
  and taking the limit $k\to\infty$,
  we obtain that 
  \begin{equation*}
    \left\| \frac{\nabla u}{|x|}\right\|_{L^2_\delta} 
    \leq \frac{1}{2} \|\Delta \tilde u\|_{L^2(\Omega)} 
    = \frac12\|\Delta u\|_{L^2(\Omega)},
  \end{equation*}
  is valid for each $\delta\in(0,R)$, 
  which implies \eqref{eq: modified Hardy-Rellich inequality n=4} by Fatou lemma used in the limit $\delta\searrow0$.
\end{proof}

\section{Global solvability. Proof of Theorem \ref{thm: global solvability}}
\label{sec: global solvability}

In this section, we are devoted to the proof of Theorem~\ref{thm: global solvability}.
We denote $C(p_1, p_2, \cdots, p_i)$ as a generic constant, 
dependent on parameters $p_1, p_2, \cdots, p_i$ only, 
which may vary from line to line.


We need the following elementary regularity estimates for the last two equations involving the Neumann heat semigroup. See e.g., \cite[Lemma~1.3]{Winkler2010b} and \cite[Lemma~2.4]{Winkler2011}.

\begin{lemma}
  \label{eq: Lp-Lq estimate of v, w}
  Let $(u,v,w)$ be a classical solution of \eqref{sys: ks isp ppp} given by Lemma~\ref{le: local existence and uniqueness}.
  Let $r\geq1$, $q\in[1,\frac{nr}{(n-2r)_+})$ and $s\in[1,\frac{nr}{(n-r)_+})$ 
  with the notation $(\varsigma)_+ := \max\{0, \varsigma\}$.
  Then there exist $C(r,q), C(r,s) > 0$ such that 
  \begin{equation*}
    \|w(\cdot, t)\|_{L^q(\Omega)} 
    \leq C(r,q)\cdot\left\{\sup_{\tau\in(0,t)}\|u(\cdot, \tau)\|_{L^r(\Omega)} + 1\right\}, 
    \quad \text{for all } t\in(0, T_{\max}),
  \end{equation*}
  and 
  \begin{equation*}
    \|\nabla v(\cdot, t)\|_{L^s(\Omega)} 
    \leq C(r,s)\cdot\left\{\sup_{\tau\in(0,t)}\|w(\cdot, \tau)\|_{L^r(\Omega)} + 1\right\},
    \quad \text{for all } t\in(0, T_{\max}),
  \end{equation*}
  holds.
\end{lemma}

For suitably weak chemotaxis strength in the sense of \eqref{eq: beta < 2n in lemma for u lp estimae},
we obtain $L^p$ estimates for $u$ by exploiting the regularity effect of the indirect signal production mechanisms, 
without imposing any restriction on $D$, e.g. slow decay or fast grow-up \eqref{h: D and S}.

\begin{lemma}
  \label{le: estimates of int u p}
  Let $n\geq2$.
  Suppose that  $\Omega, D, S$, initial data $(u_0,v_0,w_0)$ 
  and the corresponding solution $(u,v,w)$ defined in $\overline\Omega\times[0,T_{\max})$,
  are given as in Lemma~\ref{le: local existence and uniqueness}.
  Assume that  
  \begin{equation}
    0\leq S(s) \leq K_S (1+s)^{\beta-1} s, \quad s > 0,
  \end{equation}
  holds for some $K_S > 0$ and 
  \begin{equation}
    \label{eq: beta < 2n in lemma for u lp estimae}
    \beta < \frac2n.
  \end{equation}
  Then for each $p > 1$ and $T>0$, there exists $C(p, T) > 0$ such that 
  \begin{equation}
    \label{eq: finiteness of int u p}
    \int_\Omega u^p \leq C(p, T),
  \end{equation}
  holds for all $t\in(0, \min\{T, T_{\max}\})$.
\end{lemma}

\begin{proof}
  Let $p > 2 - \beta$ and put  
  \begin{equation*}
    T_0 := \min\left\{1, \frac{T_{\max}}{2}\right\}.
  \end{equation*}
  Multiplying the first equation of \eqref{sys: ks isp ppp} by $(u+1)^{p-1}$ and integrating by parts,
  we obtain 
  \begin{equation}
    \label{eq: u Lp differential inequality}
    \begin{aligned}
      \frac{1}{p}\frac{\dd}{\dd t}\int_\Omega (u+1)^p 
      &= \int_\Omega (u+1)^{p-1} \nabla\cdot (D(u) \nabla u) - \nabla \cdot(S(u)\nabla v)\\
      &= -(p-1)\int_\Omega (u+1)^{p-2}D(u)|\nabla u|^2 \\
      &\quad + (p-1)\int_\Omega (u+1)^{p-2}S(u)\nabla u \nabla v\\
      &\leq (p-1)\int_\Omega (u+1)^{p-2}S(u)\nabla u \nabla v\\
      &= (p-1) \int_\Omega \nabla F(u) \cdot \nabla v = -(p-1) \int_\Omega F(u)\Delta v\\
      &\leq p \int_\Omega F(u)|\Delta v|, \quad t \in (0, T_{\max}),
    \end{aligned}
  \end{equation}
  where 
  \begin{equation*}
    F(s) := \int_0^s(\sigma+1)^{p-2}S(\sigma)\dd\sigma,\quad s > 0.
  \end{equation*}
  Without loss of generality, we may assume 
  \[
  \beta\in\left[\frac{2}{n+2}, \frac2n\right),
  \]
  since 
  \begin{equation*}
    F(u) \leq K_S\int_0^u (1+\sigma)^{p + \beta-3}\sigma\dd\sigma \leq K_S (u+1)^{p-1+\beta}, \quad u > 0,
  \end{equation*}
  and this upper bound is an increasing function with respect to $\beta\in\mathbb R$.
  We thereby have by Young inequality 
  \begin{equation*}
    p \int_\Omega F(u)|\Delta v|
    \leq K_Sp \left(\|u+1\|_{L^p(\Omega)}^p  + \|\Delta v\|_{L^{p/(1-\beta)}(\Omega)}^{p/(1-\beta)}\right).
  \end{equation*}
  Inserting the estimate above and integrating \eqref{eq: u Lp differential inequality} over $(T_0, t) \subset (0, T_{\max})$,  
  we obtain 
  \begin{equation}
    \label{eq: u Lp estimate}
    \begin{aligned}
    &\quad \|u(\cdot, t) + 1\|_{L^p(\Omega)}^p - \|u(\cdot, T_0) + 1\|_{L^p(\Omega)}^p \\
    &\leq K_Sp^2 \left(\int_{T_0}^t\|u(\cdot, \tau) + 1\|_{L^p(\Omega)}^p \dd \tau
    + \int_{T_0}^t \|\Delta v(\cdot, \tau)\|_{L^{p/(1-\beta)}(\Omega)}^{p/(1-\beta)}\dd \tau\right),
    \quad t \in (T_0, T_{\max}).
    \end{aligned}
  \end{equation}
  By $L^p$ estimates \cite[Lemma~2.3]{Giga1991} for abstract evolution equations 
  applied to the second equation $v_t = \Delta v - v + w$ in \eqref{sys: ks isp ppp}, 
  we may pick $C(p)>0$ such that 
  \begin{equation}
    \label{eq: estimate laplace v by w}
    \begin{aligned}
    \int_{T_0}^t \|\Delta v(\cdot, \tau)\|_{L^{p/(1-\beta)}(\Omega)}^{p/(1-\beta)}\dd \tau 
    &\leq C(p)\int_{T_0}^t \|w(\cdot, \tau)\|_{L^{p/(1-\beta)}(\Omega)}^{p/(1-\beta)}\dd \tau \\
    &\quad + C(p)\|v(\cdot,T_0)\|_{W^{2,p/(1-\beta)}(\Omega)}^{p/(1-\beta)},
    \end{aligned}
  \end{equation}
  holds for all $t \in (T_0, T_{\max})$. 
  In view of conservation of mass \eqref{eq: mass identity}, 
  applying Lemma~\ref{eq: Lp-Lq estimate of v, w} to the last equation in \eqref{sys: ks isp ppp},
  we have for each 
  \begin{equation*}
    q\in\left[1,\frac{n}{\left(n-2\right)_+}\right),
  \end{equation*}
  there exists $C_q>0$ such that $\|w\|_{L^q(\Omega)} \leq C_q$.
  Using \eqref{eq: beta < 2n in lemma for u lp estimae}, 
  we can invoke Gagliardo--Nirenberg inequality and elliptic regularity theory to find $C_{gn} > 0$ and $C_e>0$ fulfilling 
  \begin{align*}
    \|w\|_{L^{p/(1-\beta)}(\Omega)}^{p/(1-\beta)} 
    &\leq C_{gn}\|w\|_{W^{2,p}(\Omega)}^p\|w\|_{L^q(\Omega)}^{p\beta/(1-\beta)}\\
    &\leq C_q^{p\beta/(1-\beta)}C_{gn}C_e \|(-\Delta + 1)w\|_{L^p(\Omega)}^p,
  \end{align*}
  for all $t\in(0, T_{\max})$,
  where $q$ is determined by 
  \begin{equation*}
    \frac{1-\beta}{p} = (1-\beta)\left(\frac1p - \frac2n\right) + \frac{\beta}q,
  \end{equation*}
  that is 
  \begin{equation*}
    q = \frac{\beta}{1-\beta}\cdot\frac n2 \in\left[1, \frac{n}{\left(n-2\right)_+}\right).
  \end{equation*}
  Another application of parabolic $L^p$ theories to the last equation in \eqref{sys: ks isp ppp} entails
  that one can find $C(p) > 0$ such that
  \begin{equation*}
    \begin{aligned}
    \int_{T_0}^t \|w(\cdot, \tau)\|_{L^{p/(1-\beta)}(\Omega)}^{p/(1-\beta)} \dd \tau
    &\leq C_q^{p\beta/(1-\beta)}C_{gn}C_e \int_{T_0}^t \|(-\Delta + 1)w(\cdot, \tau)\|_{L^p(\Omega)}^p\dd \tau\\
    &\leq C_q^{p\beta/(1-\beta)}C_{gn}C_e C(p)\left(\int_{T_0}^t \|u(\cdot, \tau)\|_{L^p(\Omega)}^p \dd \tau + \|w(\cdot, T_0)\|_{W^{2,p}(\Omega)}^p\right),
    \end{aligned}
  \end{equation*}
  holds for all $t \in (T_0, T_{\max})$.
  It follows from \eqref{eq: estimate laplace v by w} that 
  \begin{equation}
    \label{eq: estimate laplace v by u Lp}
    \begin{aligned}
    &\int_{T_0}^t \|\Delta v(\cdot, \tau)\|_{L^{p/(1-\beta)}(\Omega)}^{p/(1-\beta)}\dd \tau 
    \leq C(p)\int_{T_0}^t \|u(\cdot, \tau) + 1\|_{L^p(\Omega)}^{p}\dd \tau \\
    &\quad  + C(p) \|v(\cdot,T_0)\|_{W^{2,p/(1-\beta)}(\Omega)}^{p/(1-\beta)} 
    + C(p) \|w(\cdot, T_0)\|_{W^{2,p}(\Omega)}^p,
    \quad t\in(T_0, T_{\max}),
    \end{aligned}
  \end{equation}
  for some $C(p)>0$.
  Inserting \eqref{eq: estimate laplace v by u Lp} into \eqref{eq: u Lp estimate} 
  one can find $C(p) > 0$ such that 
  \begin{equation*}
    \|u(\cdot, t) + 1\|_{L^p(\Omega)}^p 
    \leq C(p) \int_{T_0}^t\|u(\cdot, \tau) + 1\|_{L^p(\Omega)}^{p}\dd \tau + C(p), 
    \quad t \in (T_0, T_{\max}),
  \end{equation*}
  which implies by Gronwall inequality that 
  \begin{equation*}
    \int_{T_0}^t \|u(\cdot, \tau) + 1\|_{L^p(\Omega)}^{p}\dd \tau + 1 
    \leq e^{C(p)(t-T_0)},\quad t\in(T_0, T_{\max}),
  \end{equation*}
  and thereby 
  \begin{equation*}
    \|u(\cdot, t) + 1\|_{L^p(\Omega)}^p 
    \leq C(p)e^{C(p)(t-T_0)}, \quad t \in (T_0, T_{\max}).
  \end{equation*}
  This shows \eqref{eq: finiteness of int u p}, 
  because $u$ is continuous over $\overline{\Omega}\times[0,T_0]$.
\end{proof}

\begin{proof}
  [Proof of Theorem~\ref{thm: global solvability}]
  Applying Lemma~\ref{le: estimates of int u p} for $p=n$, we can infer from Lemma~\ref{eq: Lp-Lq estimate of v, w} that for each $T>0$ and each $q>1$, there exists $C(T, q) > 0$ such that 
  \begin{equation*}
    \|w(\cdot, t)\|_{L^q(\Omega)} 
    \leq C(T,q), \quad t\in(0, \min\{T,T_{\max}\}),
  \end{equation*} 
  and thereby 
  \begin{equation} 
    \label{eq: estimate of nabla v p}
    \|\nabla v(\cdot, t)\|_{L^q(\Omega)} 
    \leq C(T,q), \quad t\in(0, \min\{T,T_{\max}\}).
  \end{equation}
  Combining with \eqref{h: D and S}, 
  we are in a position to apply a Moser-type argument \cite[Lemma~A.1]{Tao2012} to assert that for each $T>0$, 
  there exists $C(T) > 0$ such that $\|u(\cdot, t)\|_\infty \leq C(T)$ holds for all $t\in(0, \min\{T,T_{\max}\})$, 
  as a consequence of \eqref{eq: finiteness of int u p} and \eqref{eq: estimate of nabla v p}. 
  This entails $T_{\max} = \infty$ according to the extensibility criterion stated in Lemma~\ref{le: local existence and uniqueness}.
\end{proof}

\section{Lower bounds for stationary energy}
\label{sec: lower bounds for stationary energy}

The purpose of this section is to bound the stationary energy from below.
We formulate the main results as below. 
Its proof, split into two cases, 
consists of Lemma~\ref{le: mathcal F geq -C in high dimensions} and 
Lemma~\ref{le: mathcal F geq -C in 4 dimensions}.

\begin{proposition}
\label{prop: lower bounds of stationary energy n >=4}
Let $n\geq4$ and $\Omega = B_R\subset\mathbb R^n$ with some $R>0$.
Assume that  $D$ and $S$ comply with \eqref{h: D and S},
and suppose that there exist $s_0 > 1$, $\varepsilon\in(0,1)$, $K > 0$ and $k > 0$ such that 
\begin{equation*}
  \int_{s_0}^s\frac{\sigma D(\sigma)}{S(\sigma)}\dd\sigma 
\leq 
  \begin{cases}
    \frac{Ks}{\ln s},
    & \text{if } n = 4,\\
    \frac{n-4-\varepsilon}{n} \int_{s_0}^s\int_{s_0}^\sigma\frac{D(\tau)}{S(\tau)}\dd\tau\dd \sigma + Ks, 
    & \text{if } n > 4,
  \end{cases}
\end{equation*}
Then for each  $m>0$, 
there exists a constant $C=C(\varepsilon, K, s_0, R) > 0$ such that
\begin{equation*}
  \mathcal F \geq  - Cm^2 - C,
\end{equation*}
for all radial solutions $(u,v,w)\in S_m$ of stationary problem~\eqref{sys: stationary system}.
\end{proposition}

Before addressing the proof of the proposition above, 
we emphasize that the physical dimensions $(n\leq3)$ are essentially subcritical 
in the context of lower bounds of stationary energy,
which can be obtained without the contribution of diffusion and chemotaxis.

\begin{proposition}
  \label{prop: mathcal F low bounds in physical space}
  Let $n \leq 3$ and $\Omega\subset\mathbb R^n$. 
  Then there exists a constant $C > 0$ such that for each $m>0$, 
  \begin{equation*}
    \mathcal F \geq - Cm^2,
  \end{equation*}
  for all triplet $(u,v,w)\in L^1(\Omega)\times W^{2,2}(\Omega)\times L^2(\Omega)$ of nonnegative functions with $\int_\Omega u = m$ and $\partial_\nu v =0$ on $\partial\Omega$.
\end{proposition}

\begin{proof}
  The Sobolev inequalities and elliptic regularity theories 
  provide us constants $C_s, C_e>0$ such that 
  \begin{displaymath}
  \|v\|_{L^\infty(\Omega)} 
  \leq C_s\|v\|_{W^{2,2}(\Omega)} 
  \leq C_sC_e\|-\Delta v + v\|_{L^2(\Omega)},
  \end{displaymath}
  which entails 
  \begin{displaymath}
    \mathcal F(u,v,w) 
    \geq -\int_\Omega uv + \frac{1}{2}\int_\Omega|-\Delta v + v|^2 \geq - C_s^2C_e^2m^2,
  \end{displaymath}
  by Young inequality.
\end{proof}

We begin with the following direct observation.

\begin{lemma}
  \label{le: first lower bounds of stationary energy}
If $(u,v,w)$ is a solution of the stationary problem~\eqref{sys: stationary system},
then 
\begin{equation}
  \mathcal F(u,v,w) = \int_\Omega G(u) - \frac{1}{2}\int_\Omega uv 
  = \int_\Omega G(u) - \frac12\int_\Omega  |-\Delta v + v|^2,
\end{equation}
where the functional $\mathcal F$ is defined in \eqref{eq: mathcal F}.
\end{lemma}

\begin{proof}
It follows that $\int_\Omega v_t^2 = 0$ and 
\begin{align*}
  \int_\Omega  (-\Delta v + v)(-\Delta v + v)
  &= \int_\Omega w(-\Delta v + v) = \int_\Omega (-\Delta w + w)v = \int_\Omega uv.\qedhere
\end{align*}
\end{proof}

\subsection{High dimensions (\texorpdfstring{$n > 4$}{n > 4})}

To bound the stationary energy from below in high dimensions $n>4$ that is supercritical 
in the sense of the imbedding $W^{2,2}(\Omega)\subset L^\infty(\Omega)$ employed in the proof of Proposition~\ref{prop: mathcal F low bounds in physical space}, 
we first exploit the stationary problem~\eqref{sys: stationary system} to estimate $\|\Delta v\|_2$.
This is achieved by testing a fourth-order elliptic equation by $x\cdot\nabla v$ as applied in the well-known Poho\v{z}aev identity~\cite{Pohozaev1965}. 
We remark that integration by parts, combining well-designed test functions, 
was verified as a powerful tool applicable to the classification of solutions to elliptic problems~\cite{Lin2023}.

\begin{lemma}
\label{le: estimate laplace v L2 in high dimensions}
 Let $n > 4$, $m>0$, $s_0 > 0$ and $\Omega = B_R\subset\mathbb R^n$ with some $R>0$.
 Denote 
 \begin{equation}\label{eq: H}
   H(s) = H_{s_0}(s) := 
   \begin{cases}
   \int_{s_0}^s\frac{\xi D(\xi)}{S(\xi)}\dd\xi, & s \geq s_0,\\
   0, & s < s_0.
   \end{cases}
 \end{equation}
 Suppose that $(u,v,w)$ is a radially symmetric solution of \eqref{sys: stationary system} 
 with $\int_\Omega u = m$.
 Then there exists $C = C(s_0, R) > 0$ such that
 \begin{equation}\label{eq: estimate laplace v L2}
   \frac{n-4}{2}\int_\Omega |\Delta v|^2 + (n-4)\int_\Omega |\nabla v|^2 
   \leq n\int_{\{u>s_0\}}H(u)\dd x + C m^2 + C.
 \end{equation}
\end{lemma}

\begin{proof}
  Since 
    \begin{equation}
    \label{eq: four order diff v = u}
    (-\Delta + 1)(-\Delta + 1)v = (-\Delta + 1)w = u,
  \end{equation}
  motivated by \cite{Winkler2010b}, 
  we multiply this identity by $x\cdot\nabla v$ and integrate over $\Omega$ to get that 
  \begin{equation}
  \label{eq: ux grad v}
    -\int_\Omega u x\cdot\nabla v + \int_\Omega v x\cdot\nabla v
    = -\int_\Omega \Delta^2v x\cdot\nabla v
    + 2\int_\Omega \Delta v  x\cdot\nabla v.
  \end{equation}
  Two integrations by parts applied to the first term on the right in \eqref{eq: ux grad v} yield 
  \begin{align*}
    -\int_\Omega \Delta^2vx\cdot\nabla v 
    & = \int_\Omega \nabla\Delta v \cdot \nabla (x\cdot\nabla v)\\
    & = \omega_n\int_0^R r^{n-1}\left(v_{rr} + \frac{n-1}{r}v_r\right)_r \cdot (v_r + rv_{rr})\dd r\\
    & = \omega_n\left.\left(v_{rr} + \frac{n-1}{r}v_r\right)(v_rr^{n-1}+v_{rr}r^n)\right|_0^R\\
    &\quad - \omega_n\int_0^R\left(v_{rr} + \frac{n-1}{r}v_r\right)
    ((n+1)v_{rr}r^{n-1} + (n-1)v_rr^{n-2} + v_{rrr}r^n)\\
    & = 
    - \omega_n\int_0^R\left(v_{rr} + \frac{n-1}{r}v_r\right)
    \left(v_{rrr}+\frac{n-1}{r}v_{rr}-\frac{n-1}{r^2}v_r\right)r^n\dd r \\
    &\quad - 2\omega_n\int_0^R\left(v_{rr} + \frac{n-1}{r}v_r\right)^2r^{n-1}\dd r + \omega_nR^n v_{rr}^2(R)
  \end{align*}    
  and 
  \begin{align*}
  &  \quad - \omega_n\int_0^R\left(v_{rr} + \frac{n-1}{r}v_r\right)
    \left(v_{rrr}+\frac{n-1}{r}v_{rr}-\frac{n-1}{r^2}v_r\right)r^n\dd r \\
  & = - \frac{\omega_n}{2}\int_0^R\left(\left(v_{rr} + \frac{n-1}{r}v_r\right)^2\right)_rr^n\dd r\\
  & = - \frac{\omega_nR^nv_{rr}^2(R)}{2}
   + \frac{\omega_nn}{2}\int_0^R\left(v_{rr} + \frac{n-1}{r}v_r\right)^2r^{n-1}\dd r,
  \end{align*}
  because $x\cdot\nabla v = rv_r = 0$ on $\partial\Omega$.
  Therefore, we get 
  \begin{equation}\label{eq: calculate laplace2 v x grad v}
    -\int_\Omega \Delta^2vx\cdot\nabla v 
    = \frac{n-4}{2}\int_\Omega |\Delta v|^2 + \frac{\omega_nR^nv_{rr}^2(R)}{2}.
  \end{equation}
  We calculate the second term on the right side of \eqref{eq: ux grad v} by integration by parts to see that 
  \begin{equation}
  \begin{aligned}
    \int_\Omega \Delta v  x\cdot\nabla v
    &= - \int_\Omega \nabla v \cdot \nabla (x\cdot\nabla v)
    = - \omega_n\int_0^Rv_r(v_r + rv_{rr})r^{n-1}\dd r\\
    &= -\omega_n\int_0^Rv_r^2r^{n-1}\dd r 
    - \frac{\omega_n}2\int_0^R\left(v_r^2\right)_rr^n\dd r\\
    &= \frac{(n-2)\omega_n}{2}\int_0^Rv_r^2r^{n-1}\dd r = \frac{n-2}{2}\int_\Omega |\nabla v|^2.
  \end{aligned}
  \end{equation}
  As an Ehrling--type inequality enables one to find $C_R > 0$ such that 
  \begin{equation*}
    \|v\|_{L^2(\Omega)} 
    \leq \frac{1}{4R} \|\nabla v\|_{L^2(\Omega)} + C_R\|v\|_{L^1(\Omega)},
  \end{equation*}
  holds for all $v\in W^{1,2}(\Omega)$, 
  we estimate the last term on the left of \eqref{eq: ux grad v} by H\"older inequality
  and Young inequality, 
  \begin{equation}
    \begin{aligned}
    \int_\Omega v x\cdot\nabla v
    &\leq R\|\nabla v\|_{L^2(\Omega)}\|v\|_{L^2(\Omega)}
    \leq \frac{1}{4}\|\nabla v\|_{L^2(\Omega)}^2 + RC_Rm\|\nabla v\|_{L^2(\Omega)}\\
    &\leq \frac{1}{2}\|\nabla v\|_{L^2(\Omega)}^2 + R^2C_R^2m^2,
    \end{aligned}
  \end{equation}
  thanks to $m = \int_\Omega u = \int_\Omega w = \int_\Omega v$.
  In order to find an upper bound of the first term in \eqref{eq: ux grad v}, 
  we recall the definition of $H$ in \eqref{eq: H} and split the integral in question according to 
  \begin{equation}\label{eq: estimate ux grad v}
    -\int_\Omega u x\cdot\nabla v = 
    -\int_{\{u > s_0\}} u x\cdot\nabla v 
    -\int_{\{u \leq s_0\}} u x\cdot\nabla v.
  \end{equation}
  Here,  utilizing $D(u)\nabla u = S(u)\nabla v$ and integrating by parts, we obtain 
  \begin{equation}
  \label{eq: estimate ux grad v u > s0}
  \begin{aligned}
    -\int_{\{u > s_0\}} u x\cdot\nabla v 
    &= -\int_{\{u > s_0\}} x\cdot\nabla H(u)\\
    &= -\int_{\partial\Omega\cap\partial\{u>s_0\}}|x|H(u)\dd S + n\int_{\{u>s_0\}}H(u)\dd x\\
    &\leq n\int_{\{u>s_0\}}H(u)\dd x,
    \end{aligned}
  \end{equation}
  because $H(u)\equiv 0$ on $\partial\{u>s_0\}\cap\Omega$.
  Using the direct estimate 
  \begin{equation}\label{eq: estimate ux grad v u < s0}
    -\int_{\{u \leq s_0\}} u x\cdot\nabla v 
    \leq s_0R\int_{\{u\leq s_0\}}|\nabla v| 
    \leq \frac{1}{2}\|\nabla v\|^2_{L^2(\Omega)} + s_0^2R^2|\Omega|,
  \end{equation}
  we deduce \eqref{eq: estimate laplace v L2} from plugging \eqref{eq: calculate laplace2 v x grad v}--\eqref{eq: estimate ux grad v u < s0} into \eqref{eq: ux grad v}. 
\end{proof}

We now show Proposition~\ref{prop: lower bounds of stationary energy n >=4} in the case $n>4$.

\begin{lemma}
  \label{le: mathcal F geq -C in high dimensions}
Let $n > 4$, $s_0 > 1$ and $\Omega = B_R\subset\mathbb R^n$ with some $R>0$.
Suppose that 
\begin{equation}
\label{eq: H leq G high dimensions}
\int_{s_0}^s\frac{\sigma D(\sigma)}{S(\sigma)}\dd\sigma 
\leq \frac{n-4-\varepsilon}{n} \int_{s_0}^s\int_{s_0}^\sigma\frac{D(\tau)}{S(\tau)}\dd\tau\dd \sigma + Ks, 
\quad s \geq s_0,
\end{equation}
holds for some $\varepsilon\in(0,1)$ and $K>0$.
Then for each  $m>0$, 
there exists a constant $C=C(\varepsilon, K, s_0, R) > 0$ such that
\begin{equation}
\label{eq: mathcal F has lower bounds in high dimensions}
  \mathcal F \geq  - Cm^2 -C,
\end{equation}
for all radial solutions $(u,v,w)\in S_m$ of \eqref{sys: stationary system}.
\end{lemma}

\begin{proof}
  Recalling Lemma~\ref{le: first lower bounds of stationary energy},
  we shall find an upper bound for the integral
  \begin{equation*}
    \int_\Omega uv = \int_\Omega |\Delta v|^2 + 2\int_\Omega |\nabla v|^2 + \int_\Omega v^2.
  \end{equation*}
  Integration by parts and an application of H\"older inequality yield that
  \begin{equation*}
    \int_\Omega|\nabla v|^2 
    = - \int_\Omega v\Delta v  
    \leq \|v\|_{L^2(\Omega)}\|\Delta v\|_{L^2(\Omega)}.
  \end{equation*} 
  Invoking Gagliardo--Nirenberg inequality and elliptic regularity theories,
  we can find $C_{gn} > 0$ and $C_e>0$ such that 
  \begin{equation*}
    \|v\|_{L^2(\Omega)} 
    \leq C_{gn}\|v\|_{W^{2,2}(\Omega)}^\theta\|v\|_{L^1(\Omega)}^{1-\theta} 
    \leq C_{gn}C_e\|-\Delta v + v\|_{L^2(\Omega)}^\theta\|v\|_{L^1(\Omega)}^{1-\theta},
    \quad \theta = \frac{n}{n+4},
  \end{equation*}
  holds for each $v\in W^{2,2}(\Omega)$ with $\partial_\nu v = 0$ on $\partial\Omega$.
  Using triangle inequality 
  $\|-\Delta v + v\|_{L^2(\Omega)} \leq \|\Delta v\|_{L^2(\Omega)} + \|v\|_{L^2(\Omega)}$ 
  and the fact that $m = \int_\Omega u = \int_\Omega w = \int_\Omega v$,
  we can infer by Young inequality that for any $\varepsilon\in(0,1)$, there exists
  $C(\varepsilon) > 0$ such that 
  \begin{equation*}
  \label{eq: estimate uv L1}
    \int_\Omega uv \leq \frac{n-4}{n-4-\varepsilon}\int_\Omega |\Delta v|^2 + C(\varepsilon)m^2.
  \end{equation*}
  Lemma~\ref{le: first lower bounds of stationary energy} 
  and Lemma~\ref{le: estimate laplace v L2 in high dimensions} entail that 
    \begin{align*}
    \mathcal F 
    &\geq \int_\Omega G(u) 
    - \frac{1}{2}\frac{n-4}{n-4-\varepsilon}\int_\Omega |\Delta v|^2 - C(\varepsilon)m^2\\
    &\geq \int_\Omega G(u) - \frac{n}{n-4-\varepsilon}\int_{\{u > s_0\}} H(u) - C(\varepsilon, s_0, R)(m^2+1).
    \end{align*}
  As $G(s)\geq0$ for $s\leq s_0$ and by \eqref{eq: H leq G high dimensions}
  \begin{equation*}
   G(s)-\frac{n}{n-4-\varepsilon} H(s) \geqslant-\frac{n K}{n-4-\varepsilon} s,
   \quad s \geqslant s_0,
    \end{equation*}
   we can conclude that \eqref{eq: mathcal F has lower bounds in high dimensions} holds. 
\end{proof}

\subsection{Critical dimension (\texorpdfstring{$n=4$}{n=4})}

As \eqref{eq: estimate laplace v L2} indicates, 
$L^2$--estimate of $\Delta v$ for $n=4$ is critical, 
inspired by \cite{Winkler2010}, we instead use $\xi x\cdot\nabla v$ as the test function.

\begin{lemma}
  \label{le: xi' laplace v2 n = 4}
  Let $n=4$ and $\Omega = B_R\subset\mathbb R^n$ with some $R>0$.
  Then for any nonnegative and nonincreasing function $\xi(|x|)\in C^{\infty}([0,R])$ satisfying $\xi'(0) = 0$ and $\xi(R) = 0$,
  the inequality 
  \begin{equation}\label{eq: estimate laplace v 2 with xi}
  \begin{aligned} 
      - \frac32\int_\Omega \xi'|x||\Delta v|^2        
      &+ 2\int_\Omega \xi|\nabla v|^2 
      - \int_\Omega \xi'|x||\nabla v|^2\\
      &\leq 4\int_{\{u>s_0\}}\xi H(u) 
  + \int_\Omega (v + s_0)\xi |x||\nabla v| 
  + \frac12 \int_\Omega J(|x|)|\nabla v|^2,
  \end{aligned}
  \end{equation}
  holds for each radially symmetric solution $(u,v,w)$ of \eqref{sys: stationary system},
  where $H$ is defined in \eqref{eq: H} and 
  \begin{equation*}
    J(r) := -\xi'''r + 3\xi'' - \frac{3\xi'}r, \quad r > 0.
  \end{equation*}
\end{lemma}

\begin{proof}
  Motivated by~\cite{Winkler2010}, 
  we multiply \eqref{eq: four order diff v = u} by $\xi x\cdot\nabla v$ and integrate over $\Omega$ to get that
  \begin{equation}
  \label{eq: 4d u xi x grad v}
    -\int_\Omega u\xi x\cdot\nabla v + \int_\Omega v\xi x\cdot\nabla v
    = -\int_\Omega \Delta^2v\xi x\cdot\nabla v
    + 2\int_\Omega \Delta v \xi x\cdot\nabla v.
  \end{equation}
  Since $\xi(R) = x\cdot\nabla v = 0$ on $\partial\Omega$, two integrations by parts on the first term of the right yield
  \begin{equation}
    \label{eq: inte by parts twice xi x}
    - \int_\Omega \Delta^2v\xi x\cdot\nabla v
    = - \int_\Omega \Delta v\Delta(\xi x\cdot\nabla v).
  \end{equation}
  Using 
  \begin{align*}
    \Delta(\xi x\cdot\nabla v) &= \left(\partial_{rr} + \frac{n-1}{r}\partial_r\right)(\xi rv_r)\\
    &= \frac{n-1}{r}(\xi'rv_r + \xi v_r + \xi r v_{rr})  + \xi''rv_r + \xi'v_r + \xi'rv_{rr} \\
    &\quad + \xi'v_r + \xi v_{rr} + \xi'rv_{rr} + \xi v_{rr} + \xi rv_{rrr}\\
    &= \xi r\left(v_{rrr} + \frac{n-1}{r}v_{rr} - \frac{n-1}{r^2}v_r\right)\\
    &\quad + 2\xi\left(v_{rr} + \frac{n-1}{r}v_r\right) + 2\xi'r\left(v_{rr} + \frac{n-1}{r}v_r\right)\\
    &\quad + \xi''rv_r + (3-n)\xi'v_r\\
    &= \xi r(\Delta v)_r + 2\xi\Delta v + 2\xi' r\Delta v + \xi''rv_r + (3-n)\xi'v_r,
  \end{align*}
  and noting $\xi(R) = 0$,
  we calculate by integration by parts
  \begin{align*}
    - \frac{\omega_n}2\int_0^Rr^n\xi((\Delta v)^2)_r\dd r 
    &=  \frac{n\omega_n}2\int_0^Rr^{n-1}\xi|\Delta v|^2\dd r 
    +  \frac{\omega_n}2\int_0^Rr^n\xi'|\Delta v|^2\dd r \\
    &=  \frac{n}{2}\int_\Omega \xi|\Delta v|^2 
    + \frac12\int_\Omega \xi'|x||\Delta v|^2      
  \end{align*}
  and 
  \begin{align*}
    &\quad -\int_\Omega (\xi''rv_r + (3-n)\xi'v_r)\Delta v \\
    &= - \omega_n\int_0^R r^{n-1}\left(v_{rr} + \frac{n-1}{r}v_r\right)(\xi''rv_r + (3-n)\xi'v_r)\dd r \\
    &= -\omega_n\int_0^R(\xi''r^n+(3-n)\xi'r^{n-1})v_rv_{rr}\dd r \\
    &\quad - \omega_n\int_0^R\left((n-1)\xi''+\frac{(n-1)(3-n)}{r}\xi'\right)v_r^2r^{n-1}\dd r\\
    &= \frac{\omega_n}{2}\int_0^R\left(\xi'''r - (2n-5)\xi'' + \frac{(n-3)(n-1)}r\xi'\right)v_r^2 r^{n-1}\dd r
  \end{align*}
  and observe from \eqref{eq: inte by parts twice xi x}
  \begin{equation}
  \label{eq: laplace2v xi x grad v}
    \begin{aligned}
      - \int_\Omega \Delta^2v\xi x\cdot\nabla v
      &= - \frac{\omega_n}2\int_0^Rr^n\xi((\Delta v)^2)_r\dd r 
      - 2\int_\Omega \xi|\Delta v|^2 \\
      &\quad - 2\int_\Omega \xi'|x||\Delta v|^2 
      - \int_\Omega (\xi''rv_r + (3-n)\xi'v_r)\Delta v\\
      &= \frac{n-4}{2}\int_\Omega \xi|\Delta v|^2 
      - \frac32\int_\Omega \xi'|x||\Delta v|^2   
      - \frac12 \int_\Omega J_n(|x|)|\nabla v|^2,
    \end{aligned}
  \end{equation}
  with
    \begin{equation*}
    J_n(r) := -\xi'''r + (2n-5)\xi'' - \frac{(n-3)(n-1)}r\xi', \quad r > 0.
  \end{equation*}
  Integration by parts gives 
  \begin{equation*}
    -\omega_n\int_0^Rr^n\xi v_rv_{rr} \dd r
    = \frac{\omega_n}{2}\int_0^R (nr^{n-1}\xi + r^n\xi')v_r^2 \dd r,
  \end{equation*}
  and 
\begin{equation}
  \begin{aligned}
    \int_\Omega \Delta v \xi x\cdot\nabla v
    &= \omega_n\int_0^R(v_rr^{n-1})_r\xi rv_r\dd r \\
    &= - \omega_n\int_0^R v_rr^{n-1}(\xi'rv_r + \xi v_r + \xi rv_{rr})\dd r\\
    &= \frac{(n-2)\omega_n}{2}\int_0^R r^{n-1}\xi v_r^2\dd r - \frac{\omega_n}{2}\int_0^R r^n\xi'v_r^2\dd r\\
    &= \frac{n-2}{2}\int_\Omega \xi|\nabla v|^2 - \frac{1}{2}\int_\Omega \xi'|x||\nabla v|^2.
  \end{aligned}
\end{equation}
To find an upper bound for the first term in \eqref{eq: 4d u xi x grad v}, we process as in Lemma~\ref{le: estimate laplace v L2 in high dimensions} and split the integral in question according to 
\begin{equation*}
  -\int_\Omega u\xi x\cdot\nabla v 
  = - \int_{\{u>s_0\}} u\xi x\cdot\nabla v
  - \int_{\{u\leq s_0\}} u\xi x\cdot\nabla v.
\end{equation*}
Here, using $D(u)\nabla u = S(u)\nabla v$ and integrating by parts, we obtain 
\begin{align*}
  - \int_{\{u>s_0\}} u\xi x\cdot\nabla v 
  & = - \int_{\{u>s_0\}} \xi x\cdot\nabla H(u) \\
  & = n\int_{\{u>s_0\}}\xi H(u) + \int_{\{u>s_0\}}\xi'|x| H(u)\\
  & \leq n\int_{\{u>s_0\}}\xi H(u),
\end{align*}
because $H(u)\equiv0$ on $\partial\{u>s_0\}\cap\Omega$, $\xi(R)=0$ and $\xi'\leq0$.
Two direct estimates 
\begin{equation*}
  - \int_{\{u\leq s_0\}} u\xi x\cdot\nabla v
  \leq \int_\Omega s_0\xi|x||\nabla v|
\end{equation*}
and 
\begin{equation*}
  \int_\Omega v\xi x\cdot\nabla v
  \leq \int_\Omega v\xi |x||\nabla v|,
\end{equation*}
lead us to evaluate the left side of \eqref{eq: 4d u xi x grad v} by  
\begin{equation}\label{eq: estimate of left side 4d uxi}
  -\int_\Omega u\xi x\cdot\nabla v 
  + \int_\Omega v\xi x\cdot\nabla v
  \leq n\int_{\{u>s_0\}}\xi H(u) 
  + \int_\Omega (v + s_0)\xi |x||\nabla v|.
\end{equation}
Substituting \eqref{eq: laplace2v xi x grad v}--\eqref{eq: estimate of left side 4d uxi} into \eqref{eq: 4d u xi x grad v},
we arrive at 
\begin{align*}
    \frac{n-4}{2}\int_\Omega \xi|\Delta v|^2 
      - \frac32\int_\Omega \xi'|x||\Delta v|^2         
      + (n-2)\int_\Omega \xi|\nabla v|^2 
      - \int_\Omega \xi'|x||\nabla v|^2\\
      \leq n\int_{\{u>s_0\}}\xi H(u) 
  + \int_\Omega (v + s_0)\xi |x||\nabla v| 
  + \frac12 \int_\Omega J_n(|x|)|\nabla v|^2,
\end{align*}
that is \eqref{eq: estimate laplace v 2 with xi} for $n=4$.
\end{proof}

We now show Proposition~\ref{prop: lower bounds of stationary energy n >=4} in the case $n=4$.

\begin{lemma}
  \label{le: mathcal F geq -C in 4 dimensions}
Let $n=4$ and $\Omega=B_R\subset\mathbb R^n$ with some $R>0$.
Suppose 
\begin{equation}
  \label{eq: H leq Ks/ln s}
  H(s) \leq \frac{Ks}{\ln s},\quad s\geq s_0,
\end{equation}
for some $K>0$ and $s_0 > 1$.
Then for each  $m>0$, 
there exists $C(K, s_0, R) > 0$ such that
\begin{equation}
\label{eq: mathcal F has lower bounds in four dimensions}
  \mathcal F \geq  - C(K, s_0, R)m^2 - C(K, s_0, R),
\end{equation}
for all radial solutions $(u,v,w)\in S_m$ of stationary problem~\eqref{sys: stationary system}.
\end{lemma}

\begin{proof}
  Define 
  \begin{equation}\label{eq: xi}
    \xi = \xi_\eta(r) := \ln\frac{R^2+\eta}{r^2+\eta}, \quad r>0,
  \end{equation}
  with $\eta \in (0,1)$.
  Computing 
  \begin{equation*}
    \xi' = -\frac{2r}{r^2+\eta},\quad \xi''= \frac{2r^2-2\eta}{(r^2+\eta)^2},
    \quad \xi'''= \frac{12\eta r - 4r^3}{(r^2+\eta)^3},
  \end{equation*}
  we find $\xi\in C^\infty([0,R])$ is nonnegative and decreasing over $(0,R)$ with $\xi'(0) = \xi(R) = 0$.
  Therefore Lemma~\ref{le: xi' laplace v2 n = 4} is applicable to obtain 
  \begin{equation}
    \label{eq: estimate laplace v 2 with explicit xi}
    \begin{aligned} 
        3\int_\Omega \frac{|x|^2}{|x|^2+\eta}|\Delta v|^2        
        &+ 2\int_\Omega \xi|\nabla v|^2 
        + 2\int_\Omega \frac{|x|^2}{|x|^2+\eta}|\nabla v|^2\\
        &\leq 4\int_{\{u>s_0\}}\xi H(u) 
    + \int_\Omega (v + s_0)\xi |x||\nabla v| 
    + \frac12 \int_\Omega J(|x|)|\nabla v|^2.
    \end{aligned}
    \end{equation}
We calculate
    \begin{equation*}
      J(r) = - \xi'''r + 3\xi'' - \frac{3\xi'}r 
      =  \frac{16r^4}{(r^2+\eta)^3} \leq \frac{16}{r^2},
    \end{equation*}
  and estimate by Lemma~\ref{le: modified Hardy-Rellich inequality n=4} 
  \begin{equation}
    \label{eq: estimate J nabla v 2}
    \frac12 \int_\Omega J(|x|)|\nabla v|^2 \leq 8\int_\Omega\frac{|\nabla v|^2}{|x|^2} \leq 2\int_\Omega |\Delta v|^2.
  \end{equation}
  By the obvious inequality $\xi\ln\xi \geq -1/e$ for $\xi\in(0,1)$,
  we have 
  \begin{equation*}
    \xi|x| = |x|\ln\frac{R^2+\eta}{|x|^2+\eta} 
    \leq 2|x|\ln\frac{R}{|x|} \leq \frac{2R}{e},
    \quad |x| < R,
  \end{equation*}
  and infer by H\"older inequality, 
  the Ehrling-type imbedding $W^{1,2}(\Omega)\subset\subset L^2(\Omega)\subset L^1(\Omega)$,
  and Young inequality 
  that there exists $C(s_0, R) > 0$ such that 
  \begin{equation}
    \begin{aligned}
    \int_\Omega (v + s_0)\xi |x||\nabla v| 
    &\leq \frac{2R}{e}\int_\Omega (v + s_0)|\nabla v| \\
    &\leq \frac{2R}{e}\|\nabla v\|_2\|v\|_2 + \frac{2Rs_0}{e\sqrt{|\Omega|}}\|\nabla v\|_2\\
    &\leq \int_\Omega|\nabla v|^2 + C(s_0, R)(m^2 + 1).
    \end{aligned}
  \end{equation}
  Using Young-Fenchel inequality 
  \begin{equation*}
    a b \leq \frac1{\delta e}e^{\delta a} + \frac{1}{\delta}b\ln b,
  \end{equation*}
  valid for positive $a$, $b$ and $\delta$,
  we pick any $\delta\in(0,2)$ and thus estimate 
  \begin{equation}
    \label{eq: estimate xi Hu}
    \begin{aligned}
      4\int_{\{u>s_0\}}\xi H(u) 
      &= 4\int_{\{u>s_0\}}\ln\frac{R^2+\eta}{|x|^2+\eta} H(u)\\
      &\leq \frac{4}{\delta e}\int_\Omega\left(\frac{R^2+\eta}{|x|^2+\eta}\right)^\delta 
      + \frac{4}{\delta}\int_{\{u > s_0\}} H(u) \ln H(u)\\
      &\leq \frac{4R^{2\delta}}{\delta e}\int_\Omega |x|^{-2\delta} 
      + \frac{4}{\delta}\int_{\{u > s_0\}} H(u) \ln H(u), 
    \end{aligned}
  \end{equation}
  where the first integral on the right is finite since $\delta\in(0,2)$.
  \eqref{eq: H leq Ks/ln s} entails 
  \begin{equation*}
    H(s)\ln H(s) 
    \leq \frac{Ks}{\ln s}\ln \frac{Ks}{\ln s} = \frac{Ks}{\ln s}(\ln s + \ln K - \ln\ln s_0) \leq K(1+C)s,
  \end{equation*}
  for all $s > s_0$ with $C := \max\{0, (\ln K - \ln\ln s_0)/\ln s_0\}$.
  Combining with \eqref{eq: estimate J nabla v 2}--\eqref{eq: estimate xi Hu}, 
  \eqref{eq: estimate laplace v 2 with explicit xi} therefore implies 
  \begin{equation}
    \begin{aligned}
      &\quad 3\int_\Omega \frac{|x|^2}{|x|^2+\eta}|\Delta v|^2        
        + 2\int_\Omega \frac{|x|^2}{|x|^2+\eta}|\nabla v|^2\\
      &\leq \int_\Omega|\nabla v|^2 + 2\int_\Omega |\Delta v|^2 
      + C(s_0, R, K)(m^2 +1),
    \end{aligned}
  \end{equation}
  for some $C(s_0, R, K) > 0$.
  In the limit $\eta\searrow0$, Fatou lemma thus yields 
  \begin{equation*}
    \int_\Omega |\Delta v|^2        
        + \int_\Omega |\nabla v|^2
    \leq C(s_0, R, K)(m^2+1),
  \end{equation*}
  which is sufficient for boundedness of $\int_\Omega uv$ as done in Lemma~\ref{le: mathcal F geq -C in high dimensions}, 
  and thereby permits us to validate \eqref{eq: mathcal F has lower bounds in four dimensions}.
\end{proof}

\section{Initial data with large negative energy. Proof of Theorem~\ref{thm: blowup}}
\label{sec: initial data with large negative energy}

This section is devoted to existence of initial data with arbitrarily large negative energy.
We then complete the proof of Theorem~\ref{thm: blowup}.

\begin{lemma}
  \label{le: initial data with large negative energy}
  Let $n \geq 4$, $\Omega = B_R\subset\mathbb R^n$ with some $R>0$.
  Suppose that there exist $k>0$ and $s_0>1$ such that 
  \begin{equation}
    \label{eq: conditions to admit low initial energy}
    \int_{s_0}^s\int_{s_0}^\sigma\frac{D(\tau)}{S(\tau)}\dd\tau\dd \sigma 
    \leq 
    \begin{cases}
      ks(\ln s)^\theta, & \text{if } n = 4 \text{ with some } \theta\in(0,1),\\
      ks^{2-\gamma}, & \text{if } n > 4 \text{ with some } \gamma > \frac{4}{n},
    \end{cases}
  \end{equation}
  holds for all $s>s_0$.
  Then for each $m>0$ and $C>0$ one can find nonnegative $(u_0, v_0, w_0)\in (C^\infty(\overline{\Omega}))^3$ satisfying $\int_\Omega u_0 = m$ and 
  \begin{equation*}
    \mathcal F_0(u_0, v_0, w_0) < -C.
  \end{equation*}
\end{lemma}

\begin{proof}
  In the case of $n>4$, as $\gamma > 4/n$, we pick 
  \begin{equation*}
    \varrho\in\left(n- n\gamma, n-4\right)\cap(0,\infty),
  \end{equation*}
  which satisfies 
  \begin{equation}
    \label{eq: gamma criterion}
    \begin{cases}
      -\varrho < -2\varrho + n,\\
      -\varrho < -2\varrho - 2 + n,\\
      -\varrho < -2\varrho -4 + n,\\ 
      -\varrho < -(1-\gamma)n.
    \end{cases}
  \end{equation}
  Let $\phi\in C_0^\infty(\mathbb R^n)$ be a radially symmetric and nonnegative function with $\int_{\mathbb R^n}\phi = 1$, 
  which is nonincreasing with respect to $|x|$ and supported in $B_1$.
  Define for $\eta \in\left(0, R\right)$ and $\varepsilon\in(0,m)\cap(m-|\Omega|, m)$,
  \begin{equation*}
    u_\eta = \frac{m-\varepsilon}{|\Omega|} + \varepsilon \phi\left(\frac{|x|}{\eta}\right)\eta^{-n},
    \quad |x| < R.
  \end{equation*}
  and 
  \begin{equation*}
    v_\eta = \phi\left(\frac{|x|}{\eta}\right)\eta^{-\varrho},\quad w_\eta = v_\eta,\quad |x|<R.
  \end{equation*}
  Abbreviate $\|\cdot\|_p$ with $p\geq1$ as the usual norm of Lebesgue space over $B_1$.
  By direct computation, we obtain 
  \begin{equation*}
    \int_\Omega u_\eta = m, \quad   
      \int_\Omega u_\eta v_\eta = \frac{m-\varepsilon}{|\Omega|}\eta^{-\varrho + n}
      + \varepsilon\eta^{-\varrho}\|\phi\|_2^2,
  \end{equation*}
  and 
  \begin{equation*}
    \int_\Omega v_\eta^2 = \eta^{-2\varrho+n}\|\phi\|_2^2,\quad 
    \int_\Omega |\nabla v_\eta|^2 = \eta^{-2\varrho-2+n}\|\nabla\phi\|_2^2,\quad 
    \int_\Omega |\Delta v_\eta|^2 = \eta^{-2\varrho-4+n}\|\Delta\phi\|_2^2.
  \end{equation*}
  Without loss of generality, we may assume that 
  \begin{displaymath}
  \gamma\in\left(\frac{4}{n}, 1\right).
  \end{displaymath}
  We thus estimate by the basic inequality $(a+b)^c \leq 2(a^c+b^c)$, 
  valid for all $a,b>0$ and $c\in(1,2)$,
  \begin{equation*}
    \int_\Omega G(u_\eta) 
    \leq k\int_\Omega u_\eta^{2-\gamma} + \int_{\{u\leq s_0\}}G(u)
    \leq 2k\varepsilon^{2-\gamma} \eta^{-(1-\gamma)n} \|\phi\|_{2-\gamma}^{2-\gamma}
     +  C
  \end{equation*}
  with $C := 2k(m-\varepsilon)^{2-\gamma}|\Omega|^{\gamma-1} + |\Omega|A$ 
  and  $A :=\sup\{G(s): s\in((m-\varepsilon)/|\Omega|, s_0)\}$.
  According to \eqref{eq: gamma criterion}, it follows that 
  \begin{align*}
    \mathcal F_0(u_\eta, v_\eta, w_\eta) 
     & = \int_\Omega G(u_\eta) - \int_\Omega u_\eta v_\eta 
     + \frac{1}{2}\int_\Omega |\Delta v_\eta - v_\eta + w_\eta|^2 \\
     &\quad + \frac{1}{2}\int_\Omega |\Delta v_\eta|^2 + \int_\Omega |\nabla v_\eta|^2 
     + \frac{1}{2}\int_\Omega v_\eta^2\\ 
     &= \int_\Omega G(u_\eta) - \int_\Omega u_\eta v_\eta 
     + \int_\Omega |\Delta v_\eta|^2 + \int_\Omega |\nabla v_\eta|^2 
     + \frac{1}{2}\int_\Omega v_\eta^2,\\
     &\to -\infty,
  \end{align*}
  as $\eta\searrow0$, and hence the desired conclusion is valid for all sufficiently small $\eta > 0$.

  In the case of $n=4$, we leave $u_\eta$ unchanged and redefine $v_\eta$ and $w_\eta$ by 
  \begin{equation*}
    v_\eta = \psi(x)  \left(\ln\frac{R}{\eta}\right)^{-\kappa}\ln\frac{R^2+\eta^2}{|x|^2+\eta^2},
    \quad w_\eta = v_\eta, \quad |x|<R,
  \end{equation*}
  for $\eta\in(0, R/2)$,
  with 
  \begin{equation*}
    \psi := \psi_N(x) = (R^2-|x|^2)^N, \quad |x| < R,
  \end{equation*}
  for some large integer $N > 2$, where $\kappa\in(0,1)$ is small enough that $\kappa < 1 - \theta$.
  We compute  
  \begin{equation*}
    \nabla v_\eta 
    = \left(\ln\frac{R}{\eta}\right)^{-\kappa}\left(\ln\frac{R^2+\eta^2}{|x|^2+\eta^2}\nabla\psi 
    - \psi   \frac{2x}{|x|^2+\eta^2}\right),
  \end{equation*}
  and 
  \begin{equation*}
    \Delta v_\eta 
    =  \left(\ln\frac{R}{\eta}\right)^{-\kappa}\left(\Delta\psi \ln\frac{R^2+\eta^2}{|x|^2+\eta^2}
    -  \frac{4x\cdot\nabla\psi}{|x|^2+\eta^2} + \frac{4(2\eta^2 + |x|^2)\psi}{(|x|^2+\eta^2)^2}\right).
  \end{equation*}
  Writing 
  \begin{equation*}
    C_\psi := \|\psi\|_{L^\infty(\Omega)} + \|\nabla \psi\|_{L^\infty(\Omega)}
    + \|\Delta\psi\|_{L^\infty(\Omega)}, 
  \end{equation*}
  we estimate 
  \begin{equation*}
    \int_\Omega v^2_\eta \leq 4C_\psi^2(\ln 2)^{2\kappa}\int_\Omega \ln^2\frac{R}{|x|},
  \end{equation*}
  and 
  \begin{equation*}
    \int_\Omega |\nabla v_\eta|^2 
    \leq 8C_\psi^2(\ln 2)^{2\kappa}\int_\Omega \ln^2\frac{R}{|x|} 
    + 8C_\psi^2(\ln 2)^{2\kappa}\int_\Omega |x|^{-2},
  \end{equation*}
  as well as by $(a + b + c)^2\leq 3(a^2+b^2+c^2)$ valid for any $a,b,c\in\mathbb R$,
  \begin{equation*}
    \begin{aligned}
      \int_\Omega |\Delta v|^2 
      &\leq 12C_\psi^2(\ln 2)^{2\kappa}\int_\Omega\ln^2\frac{R}{|x|} 
      + 48C_\psi^2(\ln 2)^{2\kappa}\int_\Omega |x|^{-2} \\
      &\quad + 192C_\psi^2\left(\ln\frac{R}{\eta}\right)^{-2\kappa}\int_\Omega \frac{1}{(|x|^2+\eta^2)^2},
    \end{aligned}
  \end{equation*}
  where 
  \begin{equation*}
    \begin{aligned}
    &\quad \int_\Omega \frac{1}{(|x|^2+\eta^2)^2} 
    = \omega_4\int_0^{R}\frac{r^3\dd r}{(r^2+\eta^2)^2} 
    = \omega_4\int_0^{R/\eta}\frac{\xi^3\dd\xi}{(\xi^2+1)^2} \\
    & \leq \omega_4\int_0^{R/\eta}\frac{\xi\dd\xi}{(\xi^2+1)} 
    \leq \frac{\omega_4}{2}\ln\left(1+\left(\frac{R}{\eta}\right)^2\right)
    \leq \omega_4\left(\ln\frac{R}{\eta} + \frac{\ln2}{2}\right).
    \end{aligned}
  \end{equation*}
  Moreover, we estimate 
  \begin{align*}
    \int_\Omega u_\eta v_\eta 
    &\geq \varepsilon\eta^{-4}\left(\ln\frac{R}{\eta}\right)^{-\kappa} \int_\Omega \phi\left(\frac{|x|}{\eta}\right) \psi(x) 
    \ln\frac{R^2+\eta^2}{|x|^2+\eta^2}\\
    &= \varepsilon \left(\ln\frac{R}{\eta}\right)^{-\kappa} 
    \int_{B_{R/\eta}}\phi(y)\psi(y\eta)\ln\frac{R^2+\eta^2}{(|y|^2+1)\eta^2}\\
    &\geq \varepsilon \psi\left(\frac R2\right)\left(\ln\frac{R}{\eta}\right)^{-\kappa} 
    \int_{B_1}\phi(y)\ln\frac{R^2+\eta^2}{(|y|^2+1)\eta^2}\\
    &\geq \varepsilon \psi\left(\frac R2\right)\left(\ln\frac{R}{\eta}\right)^{-\kappa} 
    \left(2\ln\frac{R}{\eta} - \ln2\right),
  \end{align*}
  and by \eqref{eq: conditions to admit low initial energy},
  \begin{align*}
    \int_\Omega G(u_\eta)  
    &\leq k\int_{\{u_\eta > s_0\}} u_\eta\left(\ln u_\eta\right)^\theta + |\Omega|A   \\
    &\leq k\int_\Omega 
    \left(\frac{m}{|\Omega|} + \varepsilon \phi\left(\frac{|x|}{\eta}\right)\eta^{-4}\right) 
    \ln^\theta\left(1 + \frac{m}{|\Omega|} + \varepsilon \phi\left(\frac{|x|}{\eta}\right)\eta^{-4}\right) + |\Omega|A \\
    &\leq \frac{km}{|\Omega|}\int_\Omega \ln^\theta\left(1 + e + \frac{m}{|\Omega|} + \varepsilon \phi\left(\frac{|x|}{\eta}\right)\eta^{-4}\right)  + |\Omega|A \\
    &\quad + k\varepsilon\int_{B_{R/\eta}}\phi(y)\ln^\theta\left(1+\frac{m}{|\Omega|}\phi(y)\eta^{-4}\right) \\
    &\leq kme + \frac{km^2}{|\Omega|} + \frac{\varepsilon km}{|\Omega|} + |\Omega|A 
    + k\varepsilon \ln^\theta\left(1+\frac{m\phi(0)}{|\Omega|\eta^4}\right).
  \end{align*}
  It follows that 
  \begin{equation*}
    \mathcal{F}_0(u_\eta, v_\eta, w_\eta) 
    \leq - \frac{1}{C}\left(\ln\frac{R}{\eta}\right)^{1-\kappa} 
      + C\left(1 + \left(\ln\frac{R}{\eta}\right)^{1-2\kappa}\right) 
      + C\left(\ln\frac{R}{\eta}\right)^{\theta},
  \end{equation*}
  for all $\eta\in(0,R/2)$ with some large constant $C > 0$. 
  Our choice of $\kappa\in(0,1-\theta)$ allows us to infer that $\mathcal F(u_\eta, v_\eta, w_\eta) \to -\infty$ as $\eta\searrow0$ and conclude as before.
\end{proof}

We are in a position to show Theorem~\ref{thm: blowup}.

\begin{proof}
  [Proof of Theorem~\ref{thm: blowup}]
  For given $m>0$, 
  Proposition~\ref{prop: lower bounds of stationary energy n >=4} 
  provides us $C(\varepsilon, K, s_0, R) > 0$ such that for each stationary solution $(u_\infty,v_\infty,w_\infty)\in S_m$ of \eqref{sys: stationary system}, 
  it holds that 
  \begin{displaymath}
  \mathcal{F}(u_\infty,v_\infty,w_\infty) \geq -C(\varepsilon, K, s_0, R)m^2-C(\varepsilon, K, s_0, R).
  \end{displaymath}
  While Lemma~\ref{le: initial data with large negative energy} enables us to pick 
  nonnegative $(u_0, v_0, w_0)\in (C^\infty(\overline{\Omega}))^3$ satisfying $\int_\Omega u_0 = m$ and 
  \begin{equation*}
    \mathcal F(u_0, v_0, w_0) < -C(\varepsilon, K, s_0, R)m^2-C(\varepsilon, K, s_0, R).
  \end{equation*}
  We assume by absurdum that the corresponding solution $(u,v,w)$ exists globally and remains bounded.
  Then Lemma~\ref{le: link initial energy to stationary energy} entails that 
  there exists $(u_\infty,v_\infty,w_\infty)\in S_m$ such that 
  \begin{equation*}
    \mathcal{F}(u_\infty,v_\infty,w_\infty) - \mathcal{F}_0(u_0,v_0,w_0) 
    = - \int_0^\infty \mathcal{D}(t)\dd t \leq 0,
  \end{equation*}
  which is a contradiction.
\end{proof}

\begin{proof}
  [Proof of Corollary~\ref{coro: infinite-time blowup}]
  It suffices to validate the conditions~\eqref{h: conditions for low bounds of stationary energy} and \eqref{h: conditions for large negative initial energy} for the prototype choices of $D$ and $S$,
  which is easy to check by direct computation. 
\end{proof}

\section*{Acknowledgments}

Supported in part by National Natural Science Foundation of China (No. 12271092, No. 11671079) and the Jiangsu Provincial Scientific Research Center of Applied Mathematics (No. BK20233002).

\end{document}